\documentstyle{amsppt}

\NoBlackBoxes

\magnification=\magstep1

\define\bc{\Bbb C}
\define\br{\Bbb R}

\define\ph{\varphi}
\define\ka{\varkappa}
\define\ep{\varepsilon}

\define\car{\curvearrowright}
\define\bsl{\backslash}
\define\ovl{\overline}
\define\hlf{\frac12}
\define\inv{{-1}}

\redefine\ker{\operatorname{Ker}}
\define\tr{\operatorname{tr}}
\define\im{\operatorname{Im}}
\define\re{\operatorname{Re}}
\define\dom{\operatorname{Dom}}
\define\clos{\operatorname{clos}}
\define\dist{\operatorname{dist}}
\define\rank{\operatorname{rank}}
\define\Range{\operatorname{Range}}
\define\cC{\Cal C}
\define\spn{\operatorname{span}}

\topmatter

\title
Criteria for similarity of a dissipative integral operator to 
a normal operator
\endtitle
\rightheadtext{Similarity of a dissipative operator to a normal one}

\author S. Kupin and V. Vasyunin
\endauthor

\address
Laboratoire de Math\'ematiques Pures,\newline
Universit\'e Bordeaux 1,\newline
351, cours de la Lib\'eration,\newline
33405 Talence, France
\medskip
St. Petersburg Branch, 
V. A. Steklov Mathematical Institute,
Fontanka, 27,
191511 St. Petersburg, Russia
\endaddress

\email
kupin\@math.u-bordeaux.fr\newline
vasyunin@pdmi.ras.ru
\endemail



\abstract Let $H$ be a separable Hilbert space, $\mu$ a finite positive 
measure on $[0, 1]$. Let $\alpha$ be a measurable $L(H)$-valued function, and
$k(x,s)$ be an $L(H)$-valued positive definite kernel with 
$\tr k(x,x)\in L^1(\mu)$. 
Let, further, the values $\alpha(x)$ be selfadjoint operators. 
Sometimes we assume $\alpha(x)$ commuting with $k(x,x)$
$\mu$-almost everywhere on $[0,1]$.

We define an operator $A$ by the formula
$$
(Af)(x)=\alpha(x)f(x)+\hlf i\mu(\{x\})k(x,x)f(x)+
i\int_{[0,x)}k(x,s)f(s)d\mu(s).
$$
We are concerned with the question of similarity of $A$
to a normal operator. We obtain necessary as well as sufficient 
conditions for the similarity. The conditions turn out to be both 
necessary and sufficient in the case of a continuous measure $\mu$,
or if $\rank\im A=1$. 

Our considerations are based on function model technique. 
The key idea is to express the resolvent (and the characteristic function)
of the operator $A$ in terms of solution of a Cauchy problem and to 
analyze the spectral properties of the original operators through 
the properties of latter objects. 
This approach has been already used in~\cite{14}, \cite{16}.

A new feature of our approach is that we apply Linear Resolvent Growth 
tests~\cite{13}, \cite{15}. The test requires, in general case, the so-called 
Uniform Trace Boundedness condition  
$$
\sup_{z\in\bc_+} \im z\tr [(A^*+\ovl zI)^{-1}\im A\,(A+zI)^{-1}]<\infty.
$$
The methods allow us to study the classical problem of similarity to a 
selfadjoint operator alongside with the case of operators with non real 
discrete spectrum from basically the same point of view.

The theorems of the paper generalize to certain extent results 
of~\cite4, \cite6, \cite7, \cite{10}.

\endabstract
\endtopmatter

\head
Introduction
\endhead

One of the general goals of operator theory is to investigate spectral 
decompositions of operators  acting on a Hilbert space. For instance,
it is well known \cite{18} that an operator admits a unconditionally 
convergent spectral decomposition if and only if it is similar to a normal 
operator. Hence, to get such decompositions we need criteria, which guarantee
similarity.

Probably, the first result, obtained using the ideas we are concerned with, 
was the following theorem due to B.~Sz.-Nagy and C.~Foia\c{s}.

\proclaim{Theorem 0.1 \rm{(\cite{17})}} Let $A$ be a maximal dissipative 
operator defined on a Hilbert space and suppose that $\sigma(A)\subset \br$\,.
Then $A$ is similar to a selfadjoint operator if and only if 
$\|S_A(z)^{-1}\|\le C$ for all $z\in\bc_+$\rom, where $S_A(z)$ is the 
characteristic function of $A$ \rom(see~\thetag{1.1} for the definition\rom).
\endproclaim

The theorem can be used to obtain numerous sharp conditions of similarity for 
special classes of dissipative integral operators (cf.~\cite{10}, \cite{11}). 
We should also mention that there exist different approaches to the subject, 
see~\cite6 for a brief account.

From a certain point of view, one of the most advanced results in 
this direction was obtained by I.~Gohberg and M.~Krein.
\footnote{Added in proof. A proof of this theorem, even under more general
assumptions, appeared recently in~\cite7. A wide survey of the results
concerning similarity could be found there as well.}

\proclaim{Theorem 0.2 \rm{(\cite4)}} Let the operator 
$A^0:\;L_r^2[a,b]\to L_r^2[a,b]$ be given by the formula
$$
(A^0f)(x)=\alpha(x)f(x)+2i\int_a^x k(x,s)f(s)dm(s),
$$
where $f\in L_r^2[a,b]$\rom, $k\in C_{r\times r}([a,b]\times [a,b])$ 
is a Hermitian positive definite kernel\rom, and $\alpha$ is a real-valued 
left-continuous function. Then $A^0$ is similar to a selfadjoint operator 
if and only if the measure 
$$
\nu_c(F)=\int_{\alpha^{-1}(F)}\tr k(s,s)dm(s),
$$
is absolutely continuous and its density is uniformly bounded. Here $F$ 
denotes a measurable subset of \; $\br$ and $m$ is the Lebesgue measure.
\endproclaim

In the theorem above the indices $r$ and $r\times r$ mean that the 
corresponding spaces are composed of vector-valued or matrix-valued
functions, respectively.

The purpose of this paper is to find optimal conditions guaranteeing 
similarity of an integral operator of the same type to a normal operator. 
Namely, we consider the operator $A$ on the space $L^2(H,\mu)$ given by 
the formula
$$
(Af)(x)=\alpha(x)f(x)+i\int_0^{x+} k(x,s)f(s)d\mu(s).
\tag{0.1}
$$
Its domain will be described a bit later.
Here $\mu=\mu_c+\mu_d$ is a finite positive measure on $[0, 1]$, $\mu_c$ and 
$\mu_d$ are its continuous and discrete parts, respectively, $H$ is an 
auxiliary separable Hilbert space. Further, $\alpha$ is a $\mu$-measurable 
function on $[0, 1]$, $k$ is a positive definite kernel, and the values of 
both functions are operators on $H$. We suppose also that the values of 
$\alpha$ are selfadjoint operators on $H$ defined on a dense domain 
$\dom(\alpha(x))$.

The symbol $\int_0^{x+}$ stands for
$$
\int_0^{x+} f(t)d\mu(t)=\int_{[0,x)} f(s)d\mu(s)+\hlf\mu(\{x\})f(x).
$$

Throughout the paper we assume that $k(x,x)$ are trace class operators  
for $\mu$@-almost every $x$ and, moreover, $\tr k(x,x)\in L^1(\mu)$
(a correct meaning of the restriction of the kernel $k$ to the diagonal
will be given later in subsection~2.1).
Then the operator $A$ is well defined on the following domain
$$
\dom(A)=\left\{f\in L^2(H,\mu):\,f\in\dom(\alpha(x)\;\mu-\text{a.e., }
\int_0^1\|\alpha(x)f(x)\|_H^2d\mu(x)<\infty\right\}.
$$
Its adjoint has the same domain and is given by the formula
$$
(A^*f)(x)=\alpha(x)f(x)-i\int_{x-}^1 k(x,t)f(t)d\mu(t), 
$$
where, by definition,
$$
\int_{x-}^{1} f(s)d\mu(s)=\int_{(x,1]} f(s)d\mu(s)+\hlf\mu(\{x\})f(x).
$$
The assumptions yield that the imaginary part of $A$
$$
(2\im Af)(x)=\int_0^1k(x,s)f(s)d\mu(s)
\tag{0.2}
$$
is a trace class operator with $\tr(2\im A)=\int^1_0\tr k(x,x)\,d\mu(x)$.

We obtain necessary as well as sufficient conditions for similarity of 
$A$ to a normal operator. The conditions turn out to be both necessary and 
sufficient 
either in the case when rank of $\im A$ equals one, or in the 
case of a continuous measure $\mu$.

The proofs involve the so-called Linear Resolvent Growth tests 
(cf.~\cite{13}, \cite{15}). These tests may be viewed as strengthenings of 
Theorem~0.1. We give more details on the machinery in the next section.

The results of present work can be considered as an amplification of 
Theorem~0.2. We deal with an arbitrary measure $\mu$, we do not require 
$H$ to be of finite dimension, and we weaken the assumptions on continuity 
of $\alpha$ and $k$. Furthermore, the function $\alpha$ is 
operator-valued in our setting, and the only essential restriction on 
it is the commutativity property with respect to the kernel $k$. 
This restriction can be hardly considered as a natural one.
Without the restriction, however, we may have some new effects connected
with the non trivial discrete part of the measure $\mu$.
These effects are subject to a further investigation.

The considerations carried out in Section~3 lead to a formula
for the characteristic function of $A$. They are inspired by~\cite{16} 
and have some points in common with those in~\cite{14}.

The paper is organized as follows. In Section~1 we cite and discuss 
the theorems we rely on. We formulate the main results of the paper in 
Section~2. The proofs and an example are presented in Section~3.

We finish the introduction with some notation. If $H$ is a Hilbert space, 
we denote the space of all linear bounded operators on $H$ by $L(H)$. 
We recall that a densely defined operator $A$ with domain $\dom(A)=\dom(A^*)$ 
acting on a Hilbert space $H$, is called dissipative if $\im A=
\frac{A-A^*}{2i}\ge 0$. The spectrum of $A$ is denoted by $\sigma(A)$ and 
$\sigma_p(A)$ stands for its point part. We say that a dissipative operator 
$A$ is maximal if it has no dissipative extensions. Operators $A$ and $B$ 
are called similar if there exists a bounded and boundedly invertible 
operator $V$ such that $A=V^{-1}BV$. The symbol $\goth S_1$ stands, as
usual, for the ideal of operators, and $\goth S_2$ for the
ideal of Hilbert--Schmidt operators. We shall write ``a.e.'' instead of 
$\mu$@-a.e.

\head
1. Resolvent tests for dissipative operators
\endhead

\subheading{1.1. Resolvent tests} 
It is well known (see~\cite{12}, \cite{16}) that a maximal dissipative 
operator $A$ admits  a canonical representation  of the form 
$A=A_0\oplus A_1$, where $A_0$ is a selfadjoint operator and $A_1$ is a 
completely non-selfadjoint one. Since $A_0$ is itself normal, it is  
obvious that $A$ is similar to a normal operator if and only if $A_1$ is. 
The operator $A_1$ is unitarily equivalent to the Sz.-Nagy-Foia\c{s} model 
operator (see~\cite{12}, \cite{16}), constructed with the help of the 
so-called characteristic function $S_{A_1}$ of $A_1$ (see~\thetag{1.1} below). 
Formula~\thetag{1.1} shows, in particular, that the characteristic 
functions of $A$ and $A_1$ coincide. Hence, the question of similarity 
of $A$ to a normal operator is reduced to the same question for its
completely non-selfadjoint part $A_1$ given in terms of $S_A$. The 
selfadjoint part $A_0$ can be merely ``dropped''.

We have the following theorem.
\proclaim{Theorem 1.1 \rm{(\cite{15})}} Let $A$ be a maximal dissipative 
operator defined on a Hilbert space and suppose that $\sigma(A)\neq\ovl\bc_+$. 
Then $A$ is similar to a normal operator if the following conditions hold
$$
\align
\text{\rm i)}\ &\cC_1(A)=\sup_{z\in\bc_+\bsl\sigma(A)} 
\|(A-zI)^{-1}\|\cdot\dist(z,\sigma(A))<\infty,
\tag{LRG}\\
\text{\rm ii)}\ &\cC_2(A)=\sup_{z\in\bc_+} 4\im z\cdot
\tr\big[(A^*-zI)^{-1}\im A\;(A-\ovl zI)^{-1}\big]<\infty.
\tag{UTB}
\endalign
$$
\endproclaim

We write (LRG) and (UTB) instead of Linear Resolvent Growth and 
Uniform Trace Boundedness for the sake of brevity.
We see that for $z\in\bc_+$ 
$$
I-b_z(A)^*b_z(A)=4\im z\cdot(A^*-zI)^{-1}\im A\;(A-\ovl zI)^{-1}\ge 0,
$$
where $b_z(w)=(w-z)/(w-\ovl z)$ is an elementary Blaschke factor in $\bc_+$. 
Consequently,
$$
\rank (I-b_z(A)^*b_z(A))=\rank \im A.
$$
Hence, if $n=\rank \im A<\infty$, then $\cC_2(A)\le n<\infty$, and the (UTB) 
condition is satisfied automatically. The latter case was studied alongside 
with (LRG) tests for operators with Dini-smooth spectral 
sets (and hence for dissipative ones, too), in \cite{13}.
The previous theorem can be viewed as a generalization of a 
result from that work. 

Note that Theorem~1.1 was originally proved for contractions. The version 
we present here can be easily obtained from the main theorem of~\cite{15} by 
applying the Cayley transform
$$
T=(A-iyI)(A+iyI)^{-1},\quad A=iy(I+T)(I-T)^{-1},\quad y>0.
$$

Further, Theorem~1.1 admits a conformally invariant transcription. 
We need to define the characteristic function of $A$  to formulate 
the result (see \cite{12}, \cite{16})
$$
S_A(z)=I+i(2\im A)^\hlf(A^*-zI)^{-1}(2\im A)^\hlf|_{\Range\im A}.
\tag{1.1}
$$
We recall that $\|S_A(z)\|\le 1$ for $z\in\bc_+$.

The following theorem is implicitly contained in \cite{15}.

\proclaim{Theorem 1.2} Let $A$ be a maximal dissipative operator on a 
Hilbert space and suppose that $\sigma(A)\ne\ovl\bc_+$. Then $A$ 
is similar to a normal operator if the following conditions hold
$$
\align
\text{\rm i)}\ &\cC_3(A)=\sup_{z\in\bc_+\bsl\sigma(A)}\{\|S_A(z)^{-1}\|
\cdot\inf_{\lambda\in\sigma_p(A)\cap\bc_+}|b_\lambda(z)|\}<\infty,
\tag{1.2}\\
\text{\rm ii)}\ &\cC_2(A)=\sup_{z\in\bc_+}\tr(I-S_A(z)^*S_A(z))<\infty.
\tag{1.3}
\endalign
$$
\endproclaim

If $\sigma_p(A)\cap\bc_+=\emptyset$, we understand inequality~\thetag{1.2} as
$$
\cC_3(A)=\sup_{z\in\bc_+} \|S_A(z)^{-1}\|<\infty.
$$

The advantage of the following theorem is that it does not require 
the (UTB) condition in the sufficient part of the statement.

\proclaim{Theorem 1.3} Let $A$ be a maximal dissipative operator with
a trace class imaginary part. Then $A$ is similar to a normal 
operator if 
\roster
\item"i)" there is no singular inner part in the canonical factorization of 
the function $\det S_A(z)$ and the outer part $\{\det S_A(z)\}_{out}$ is 
bounded away from zero
$$
|\{\det S_A(z)\}_{out}|=|\{\det S_A(z)\}_{sing,\; out}|\ge\delta>0,
\qquad z\in\bc_+;
\tag{1.4}
$$
\item"ii)" the operator $A$ has no root vectors and 
the family of subspaces $\{\ker (A-zI)\}$, $z\in\sigma_p(A)\cap\bc_+$\rom,
forms an unconditional basis in its linear span.
\endroster
Condition {\rm ii)} is always necessary for the operator $A$ to be similar
a normal operator. Condition {\rm i)} is necessary under assumption
that $A$ has the~\thetag{UTB} property.
\endproclaim

\demo{Sketch of the proof} First, we note that $I-S_A\in\goth S_1$
for the dissipative operators with a trace class imaginary part, i.e.,
the determinant $\det S_A$ is well defined.

We represent $S_A$ as a product of two factors: first one is
a Blaschke--Potapov product, the second one is invertible in $\bc_+$. This
factorization single out two invariant subspaces, say $H_1$ and $H_2$,
such that $\sigma(A|H_1)\subset\br$ and 
$A|H_2$ has the point spectrum only.
It is clear that $A$ is similar to a normal operator if and only if
\roster
\item"a)" $A|H_1$ is similar to a selfadjoint operator;
\item"b)" $A|H_2$ has no root spaces and the family of subspaces 
$\{\ker (A|H_2-zI)\}$, $z\in\sigma_p(A|H_2)$, forms an unconditional 
basis in $H_2$;
\item"c)" the angle between $H_1$ and $H_2$ is positive and their
direct sum is the whole $H$.
\endroster

By the Sz.-Nagy--Foia\c{s} result (Theorem~0.1), the claim a) is 
equivalent to say that $S_A$ has no singular inner factor
and $\|\{S_A(z)\}_{out}^{-1}\|\le C$ for all $z\in\bc_+$.
Recall the following inequality
$$
\exp\left\{-\|T^\inv\|\tr(I-T)\right\}\le\det T\le\frac1{\|T^\inv\|}
$$
valid for any positive contraction $T$ (sf.~\cite{15}, Lemma~4.1). 
We see that if determinant of $\{S_A(z)\}_{out}$ is bounded away from zero, 
then operator $\{S_A(z)\}_{out}^{-1}$ is uniformly bounded in $\bc_+$.
The latter condition is also sufficient to bound $\det \{S_A(z)\}_{out}$ 
away from zero, if the (UTB) condition is fulfilled.

The claim b) is the same as ii). Finally,
using~\cite2, we easily obtain that i) implies c).
\qed
\enddemo

\subheading{1.2. $N$-Carleson sets and Carleson measures}
In this subsection we recall some facts about geometry of point subsets 
in $\bc_+$.

For a point $z$, $z\in\bc_+$, and a number $\delta$, $0<\delta<1$, we put
$$
B_\delta(z)=\{w\in\bc_+:\; |b_z(w)|\le\delta\}.
$$
Let $\Lambda=\{z_k\}$ be a point sequence in the upper half-plane $\bc_+$. 
We say that the sequence $\Lambda$ is sparse, if
$$
\inf\{|b_{z_k}(z_j)|: z_k,z_j\in\Lambda,\ k\ne j\}>0,
$$
or, equivalently, if there exists a number $\ep>0$ such that $B_\ep(z)\cap 
B_\ep(w)=\emptyset$ for $z,w\in\Lambda$ and $z\ne w$.

The sequence is called Carleson, if
$$
\inf_{z\in\Lambda}\prod_{w\in\Lambda\bsl\{z\}}|b_w(z)|\ge\delta_0>0.
\tag{1.5}
$$
It can be shown (see \cite8, ch.\;9) that the latter is equivalent to say 
that there exists a number $c=c(\delta_0)>0$ with the property
$$
\prod_{w\in\Lambda}|b_w(z)|\ge c\inf_{w\in\Lambda}|b_w(z)|
\tag{1.6}
$$
for all $z\in\bc_+$.  
If $\Lambda=\cup^N_{k=1}\Lambda^k$, where $\Lambda^k$ are sparse (Carleson),
we say that $\Lambda$ is $N$@-sparse ($N$@-Carleson).
Further, we say that a measure $\sigma$ on $\bc_+$ is a Carleson measure, 
if 
$$
\sigma(Q)\le Ch
$$
for a constant $C$ and all squares $Q=[x-h, x+h]\times i[0,2h]$,
$x\in\br$.

Detailed discussion of notions we mentioned can be found
in~\cite8, ch.~7. We only state an equivalence we
use in the sequel.

\proclaim{Theorem 1.4}
Let $\Lambda\subset\bc_+$. The following assertions are equivalent.
\roster
\item"i)" The measure $\sigma=\sum_k\im z_k\,\delta_{z_k}$ is Carleson\rom;
\item"ii)"
$\displaystyle
\sup_{z\in\bc_+}\sum_k\frac{\im z\im z_k}{|z-\ovl{z_k}|^2}<\infty,
$
\newline
the symbol $\delta_{z_k}$ stands for Dirac's delta measure at the 
point $z_k$\rom; 
\item"iii)" $\Lambda$ is $N$@-Carleson for some $N$.
\endroster
\endproclaim

\subheading{Remark 1.5}. Note that $\Lambda$ is $N$@-Carleson for a given $N$
if and only if the corresponding measure is Carleson and $\Lambda$ is
$N$@-sparse with the same N.

\medskip
The following result is stated as Corollary~3.3 of~\cite{15}.

\proclaim{Theorem 1.6} Under the \thetag{LRG} and the \thetag{UTB} conditions
$\sigma_p(A)\cap\bc_+$ is a $N$@-Carleson set for some $N<\infty$.
\endproclaim

\head
2. The main results
\endhead

\subheading{2.1. Notation}
It follows from formula~\thetag{0.2} that the operator $A$ given 
by~\thetag{0.1} is dissipative. It will be shown in Corollary~3.11 that 
$A$ is maximal.

Following~\cite{16}, we define the functions
$$
\align
\ph: [0,1]\to[0,M],\qquad 
\ph(x)&=
\cases
\hlf\mu(\{0\}),& x=0,\\
\mu([0,x))+\hlf\mu(\{x\}),& x>0;
\endcases\\
\psi: [0,M]\to [0,1],\qquad
\psi(t)&=
\cases
\inf\{x:\mu([0,x))>t\},& t<\mu([0,1))\\
1,& t\ge\mu([0,1)),
\endcases
\endalign
$$
where we have put $M=\mu([0,1])$. Given a function $f$ on $[0,1]$, we set 
$$
f_*(t)=f(\psi(t)),\qquad t\in[0,M].
$$
For instance, $\ph_*(t)$ stands for $\ph(\psi(t))$. 
Notice that $\ph_*(t)=t$ if $\mu(\{\psi(t)\})=0$.
Observe also that $\psi(s)=const=x$ on the interval 
$(\ph(x-0),\ph(x+0))$ and, consequently, $v_*(s)=v(x)$ on 
$(\ph(x-0),\ph(x+0))$ for every $v$ defined on $[0,1]$.

For the sake of brevity we often write $\mu_x$ instead of 
$\mu(\{x\})$ whenever $\mu(\{x\})>0$.
 
We mention that if $f\in L^1(\mu)$, a change of 
variables yields
$$
\int_0^{\ph(x)}f_*(t)dt=\int_0^{x+}f(s)d\mu(s).
\tag{2.1}
$$

We introduce some more objects. Let $E$ be an auxiliary Hilbert space
of dimension $\dim\Range\im A$. We take an operator $c$ from $E$ 
to $L^2(H,\mu)$ with the property $cc^*=2\im A$. 
Since $\im A\in\goth S_1$, $c$ is a Hilbert--Schmidt operator and it 
defines an operator-valued function $\mu$-almost everywhere (we shall 
use for the function the same notation $c$). Note that $(ch)(x)=c(x)h$, 
$h\in E$. The values $c(x)$ are Hilbert--Schmidt 
operators from $E$ to $H$ and the kernel $k$ in the 
definition~\thetag{0.1} of the operator $A$ can be written as 
$$
k(x,s)=c(x)c(s)^*.
$$
To check the above assertion, we choose an orthonormal basis, say, $\{e_j\}$
in $E$ and put $c_j=ce_j$, $c_j\in L^2(H,\mu)$. Since $c\in\goth S_2$,
$$
\|c\|^2_{\goth S_2}=\sum\|ce_j\|^2_{L^2(H,\mu)}
=\sum\int^1_0\|c_j(x)\|^2_H\,d\mu(x)<\infty,
$$
and therefore, the formula
$$
c(x)\bigl(\sum y_je_j\bigr)=\sum y_jc_j(x)
$$
defines a bounded operator from $E$ to $H$ for almost all $x$. 
Furthermore,
$c(x)\in\goth S_2$ a.e. Indeed,
$$
\|c(x)\|^2_{\goth S_2}=\sum\|c(x)e_j\|^2_H=\sum\|c_j(x)\|^2_H<\infty
\quad\text{for a.e. } x,
$$
because, as we have seen, the function $x\mapsto\|c(x)\|_{\goth S_2}$
belongs to $L^2(\mu)$. Since the adjoint operator $c^*$ is given by the formula
$$
c^*f=\int^1_0c(x)^*f(x)d\mu(x), \qquad f\in L^2(H,\mu),
$$
we can easily check the relation
$$
k(x,s)=c(x)c(s)^*=\sum(\cdot,c_j(s))_Hc_j(x)
$$
by means of~\thetag{0.2}:
$$
\align
\int^1_0\int^1_0&\big(k(x,s)f(s),g(x)\big)_H\,d\mu(s)d\mu(x)
=(2\im Af,g)_{L^2(H,\mu)}=(c^*f,c^*g)_E\\
&=\int^1_0\int^1_0\big(c(s)^*f(s),c(x)^*g(x)\big)_E\,d\mu(s)d\mu(x)\\
&=\int^1_0\int^1_0\big(c(x)c(s)^*f(s),g(x)\big)_H\,d\mu(s)d\mu(x).
\endalign
$$
So, we write $k(x,x)$ assuming $c(x)c(x)^*$.

Furthermore, using the spectral theorem for the selfadjoint operator $\alpha(x)$
in von Neumann form we can represent the space $H$ as a direct integral
$$
H=\oplus\int_\br H_x(\lambda)d\rho_x(\lambda)
$$
with respect to the scalar spectral measure $\rho_x$ of the operator 
$\alpha(x)$. Now we can define Hilbert--Schmidt operators $c(x,\lambda)$ 
acting from $E$ to $H_x(\lambda)$ in the same way as the operators $c(x)$
were defined.
We choose an arbitrary orthonormal basis $\{e_j\}$
in $E$ and put $c_j(x,\lambda)=(ce_j)(x,\lambda)$, 
$c_j(x,\lambda)\in H_x(\lambda)$. Since $c\in\goth S_2$,
$$
\|c\|^2_{\goth S_2}=\sum\|ce_j\|^2_{L^2(H,\mu)}
=\sum\int^1_0\int_\br\|c_j(x,\lambda)\|^2_{H_x(\lambda)}
d\rho_x(\lambda)\,d\mu(x)<\infty,
$$
and therefore, the formula
$$
c(x,\lambda)\bigl(\sum h_je_j\bigr)=\sum h_jc_j(x,\lambda)
$$
defines a Hilbert--Schmidt operator from $E$ to $H_x(\lambda)$ for 
$\rho_x$@-almost all $\lambda$ and $\mu$@-almost all $x$. 
For the adjoint operator we have the formulas:
$$
\align
c(x)^*f(x)&=\int_\br c(x,\lambda)^*f(x,\lambda)\,d\rho_x(\lambda),\\
c^*f&=\int^1_0\int_\br c(x,\lambda)^*f(x,\lambda)\,d\rho_x(\lambda)\,d\mu(x)
\endalign
$$
for arbitrary $f\in L^2(H,\mu)$.

In some theorems we assume the following commutativity relation
$$
k(x,x)[\alpha(x)-zI]^\inv=[\alpha(x)-zI]^\inv k(x,x)
\tag{2.2}
$$
to be fulfilled. This will allow us to separate the point part of the 
spectrum from its real continuous part. Under this assumption every 
eigensubspace of $k(x,x)$ is invariant under $\alpha(x)$. Since all 
these subspaces corresponding to non-zero eigenvalue of $k(x,x)$ are 
finite dimensional, we can choose there an orthogonal family of 
eigenvectors of $\alpha(x)$ that by construction are eigenvectors of 
$k(x,x)$ as well. So, we can construct a sequence of orthonormal vectors 
$e_j(x)$ such that
$$
\align
k(x,x)&=\sum\kappa_j(x)^2(\,\cdot\,,e_j(x))_He_j(x),\\
\alpha(x)&=\sum\alpha_j(x)(\,\cdot\,, e_j(x))_He_j(x)+\alpha_0(x),
\tag{2.3}
\endalign
$$
where $\alpha_0(x)=\alpha(x)|_{\ker k(x,x)}$ has no importance for us
because it did not touch the completely non-selfadjoint part and is
included in the selfadjoint part of our operator $A$.

\subheading{2.2. The characteristic function of $A$} We shall calculate the 
characteristic function $S_A$ as a solution of a certain Cauchy problem. 
Recall that the characteristic function is defined up to constant unitary 
operators from the left and from the right. It will be convenient 
to consider $S_A$ on the auxiliary Hilbert space $E$ rather than on 
$\clos\Range\im A$. To do this, we map $E$ onto $\clos\Range\im A$ 
with the help of the 
isometry $V_c$ from the polar decomposition $c=(2\im A)^\hlf V_c$.
Now,  we can calculate $S_A$ as
$$
S_A(z)=I+ic^*(A-zI)^{-1}c
\tag{2.4}
$$

Define a $L(E)$@-valued function $G$ as the solution of the following 
Cauchy problem 
$$
\cases
G(t,z)'&=\;\;c_*(t)^*\Omega(t,z)c_*(t)G(t,z)\\
G(M,z)&=\;\;I,
\endcases
\tag{2.5}
$$
where 
$$
\Omega(t,z)=\Bigl[\big(t-\ph_*(t)\big)k_*(t,t)+i(\alpha_*(t)-zI)\Bigr]^{-1}
\tag{2.6}
$$
and $k_*(t,t)=c_*(t)c_*(t)^*=k(\psi(t),\psi(t))$.

\proclaim{Theorem 2.1} Cauchy problem~\thetag{2.5} is solvable for
$$
\im z\ge1+\hlf\sup\{\mu_x\|k(x,x)\|:\,\mu_x>0\}
\tag{2.7}
$$
\rom(for $\im z\ge1$ when $\mu=\mu_c$\rom) and
$$
S_A(z)=G(0,z).
$$
\endproclaim

The following corollary will be deduced from the theorem by straightforward 
calculations.

\proclaim{Corollary 2.2} $I-S_A(z)\in\goth S_1$ and under 
assumption~\thetag{2.2} we have
$$
\det S_A(z)=\ep\prod_{j,x:\,\mu_x>0} 
\left(\frac{z-z_j(x)}{z-\ovl{z_j(x)}}e^{i\phi_j(x)}\right)\cdot
\exp\left(i\int_0^1\tr\left[c(x)^*(\alpha(x)-zI)^{-1}c(x)\right]
d\mu_c(x)\right).
\tag{2.8}
$$
Here $\ep$ is a unimodular constant,
$z_j(x)=\alpha_j(x)+\hlf i\mu_x\kappa_j^2(x)$\rom, where $\{\kappa_j^2(x)\}$
and $\{\alpha_j(x)\}$ are the set of positive eigenvalues of $k(x,x)$ 
(or\rom, what is the same\rom, of $c(x)^*c(x)$) and the set of corresponding
eigenvalues of $\alpha(x)$ on a common family of eigenvectors. The real 
numbers $\phi_j(x)$ are chosen to satisfy 
$e^{i\phi_j(x)}{(i-z_j(x))}/{(i-\ovl{z_j(x)})}>0$.
\endproclaim

It follows that the non real part of the point spectrum of $A$
$\sigma_p(A)\cap\bc_+$ coincides with the set $\{z_j(x)\}_{x:\,\mu_x >0}$.

We note that the Blaschke product in~\thetag{2.8} is well defined. Indeed,
$$
\align
\sum_{j,\,x:\,\mu_x>0}\frac{\im z_j(x)}{1+|z_j(x)|^2}
&\le\sum_{j,\,x:\,\mu_x>0} \im z_j(x) 
=\hlf\sum_{j,\,x:\,\mu_x>0} \mu_x\kappa_j^2(x)\\
&=\hlf\sum_{x:\,\mu_x>0} \mu_x\tr k(x,x)
\le\hlf\int_0^1\tr k(x,x)d\mu(x)<\infty,
\tag{2.9}
\endalign
$$
since the function $x\mapsto\tr k(x,x)$ lies in $L^1(\mu)$.

\subheading{2.3. Necessary conditions} Corollary~2.2 and the results 
of Section~1 allow us to obtain the following necessary conditions for 
similarity of $A$ to a normal operator. The conditions are 
also sufficient for similarity for operators with a continuous measure $\mu$. 

\proclaim{Theorem 2.3} If the operator $A$ defined by~\thetag{0.1} and
satisfying \thetag{UTB} condition 
is similar to a normal operator, then the measure
$$
\nu_c(F)=\int_0^1\int_F
\tr c(x,\lambda)^*c(x,\lambda)\,d\rho_x(\lambda)\,d\mu_c(x),
$$ 
is absolutely continuous with respect to Lebesgue measure $m$\rom, and 
its density is uniformly bounded\rom, i.e.\rom,
$$ 
d\nu_c(s)=w(s)dm(s), \quad w\in L^\infty(m).
\tag{2.10}
$$ 
Here $F$ is an arbitrary Borel subset of\; $\br$.

Conversely\rom, if the measure $\mu$ is continuous \rom(i.e.\rom, 
$\mu=\mu_c$\rom) and the density of the measure $\nu_c$ with respect 
to Lebesgue measure is uniformly bounded\rom, then the operator $A$ is
similar to a selfadjoint operator even without \thetag{UTB} restriction.
\endproclaim

As to the discrete part of the spectrum we, of course,
can say that it satisfies the Blaschke condition~\thetag{2.9}. This 
fact has nothing in common with similarity to a normal operator, it is 
just
a consequence of the fact that the imaginary part of $A$ is a trace
class operator. We are unable to state other reasonable conditions 
necessary for similarity in terms of the measure $\mu$, the function 
$\alpha$, and the kernel $k$ without some additional assumptions. 
Nevertheless, under the (UTB) condition and assuming that $\alpha(x)$ 
commute with $k(x,x)$ we can formulate the following assertion that also 
looks like boundedness of a density of a measure.

\proclaim{Theorem 2.4} Suppose that the operator $A$ defined by 
formula~\thetag{0.1} and satisfying~\thetag{2.2} is similar to a normal 
operator and that the \thetag{UTB} condition is true. Then the family 
of measures
$$
\nu_{d,h}(F)=\sum_j\sum_{x\in\alpha_j^{-1}(F)}\eta_{4h}(\mu_x\kappa_j(x)^2)
$$
satisfies the condition
$$
\sup_{h>0}\;\sup_{x_0\in\br}\frac{\nu_{d,h}([x_0-h,x_0+h])}h<\infty, 
\tag{2.11}
$$
where
$\eta_t(x)=x\cdot\chi_{[0,t]}(x)$\rom, and $\chi_{[0,t]}(x)$ is the 
indicator of the interval $[0,t]$.
\endproclaim

The proof of  the theorem will show (see Lemma 3.13) that 
condition~\thetag{2.11} is nothing but the condition for the 
canonical measure $\sigma$ related with the point spectrum of $A$ to be 
Carleson, or, in other words, the condition for the spectrum itself to 
be an $N$-Carleson set (see Theorem~1.4). This property was already 
stated exactly in this form in Theorem~1.6. We shall see after 
studying the (UTB) condition (Corollary~2.9), that if the non-real
spectrum of $A$ is $N$@-Carleson and~\thetag{2.10} takes place, 
then the (UTB) condition is fulfilled.

Conditions \thetag{2.10} and \thetag{2.11} of the previous two theorems 
can be combined in the following manner. Define a measure  $\nu_h$
by the formula
$$
\nu_h=\nu_c+\nu_{d,h}.
$$
Now, conditions \thetag{2.10} and \thetag{2.11} hold if and only if
$$
\sup_{h>0}\;\sup_{x_0\in\br}\frac{\nu_h([x_0-h,x_0+h])}h<\infty, 
\tag{2.12}
$$

\subheading{2.4. Sufficient conditions} We begin the subsection with 
several general observations. It is known that for every $z\in\bc_+$ 
$$
\dim\ker (A^*-\ovl{z}I)=\dim\ker S_A(z)^*,\qquad
\dim\ker (A-zI)=\dim\ker S_A(z).
$$
Since, by Corollary~2.2, the determinant of $S_A$ exists, we always have
$\dim\ker S_A(z)^*=\dim\ker S_A(z)$, and, consequently,
$\dim\ker (A-zI)=\dim\ker (A^*-\ovl{z}I)$. 

\proclaim{Theorem 2.5} Suppose that the operator $A$ given 
by~\thetag{0.1} and satisfying~\thetag{2.2}
has the following properties:
\roster
\item"i)" inequality~\thetag{2.12} holds \rom(or\rom, equivalently\rom, 
conditions~\thetag{2.10} and~\thetag{2.11} are satisfied\rom)\rom;
\item"ii)" $\sigma_p(A)\cap\bc_+$ is a sparse sequence\rom;
\item"iii)" the operator $A$ does not have root vectors\rom;
\endroster
then $A$ is similar to a normal operator.
\endproclaim

In the case $\rank\im A=1$ the conditions of the above theorem
turn out to be necessary. Since condition~iii) follows from~i) and~ii),
we obtain the following criterion.

\proclaim{Theorem 2.6} Let $A$ be the operator given by~\thetag{0.1}
with $\rank\im A=1$ and satisfying~\thetag{2.2}
\rom(i.e.\rom, $c(x)$ is an eigenvector of $\alpha(x)$\rom).
Then $A$ is similar to a normal operator if and 
only if the following conditions are fulfilled:
\roster
\item"i)" condition~\thetag{2.12} holds;
\item"ii)" the sequence $\{z_x\}_{x:\, \mu_x>0}$ is sparse.
\endroster
\endproclaim

The deficiencies of Theorem~2.5, from our point of view, are as follows: 
first, condition~ii)
cannot be necessary in more or less general situation, and, second,
though condition~iii) is necessary for $A$ to be similar to a normal
operator, it is not expressed via natural terms of the problem, 
i.e., via terms of the measure $\mu$, the function $\alpha$, and the 
kernel $k$. We have to know more about geometry of eigen and root
spaces in these parameters. The following lemma is the first step in 
this direction. Assertions of this type might permit to 
strengthen Theorem~2.5 and to obtain more precise sufficient conditions 
for similarity of $A$ to a normal operator.

\proclaim{Lemma 2.7}
For every $x\in[0,1]$ with $\mu_x>0$, the function $S_A$ admits the following 
regular factorization
$$
S_A(z)=S_{x-}(z)B_x(z)S_{x+}(z),
\tag{2.13}
$$
here $B_x(z)=\left[I+\frac{i}2\mu_xc(x)^*(\alpha(x)-zI)^{-1}c(x)\right]$ 
$\left[I-\frac{i}2\mu_xc(x)^*(\alpha(x)-zI)^{-1}c(x)\right]^{-1}$\rom,
and other two factors can be expressed in 
terms of the solution $G$ of Cauchy problem~\thetag{2.5}\rom:
$S_{x+}(z)=G(\ph(x+0),z)$, $S_{x-}(z)=G(0,z)G(\ph(x-0),z)^{-1}$.
\endproclaim

For readers familiar with the notion of multiplicative integral 
(see~\cite9), we mention that the factors in~\thetag{2.13}
can be written as follows:
$$
\align
S_{x-}(z)&=\int^{\ph(x-0)\atop\car}_0 e^{-c_*(t)^*\Omega(t,z)c_*(t)dt},\\
B_x(z)&=
\int_{\ph(x-0)}^{\ph(x+0)\atop\car}e^{-c_*(t)^*\Omega(t,z)c_*(t)dt},\\
S_{x+}(z)&=\int^{M\atop\car}_{\ph(x+0)} e^{-c_*(t)^*\Omega(t,z)c_*(t)dt}.
\endalign
$$

We shall say more about some  factorization properties of $S_A$ after we
prove this lemma and introduce additional notation.

\subheading{2.5 More about (UTB)} A lemma  we  formulate in this subsection 
provides an alternative way  of calculating the constant $\cC_2(A)$. 
It is convenient in verifying the (UTB) property.

\proclaim{Lemma 2.8} Let $G$ be the solution of the Cauchy 
problem~\thetag{2.5}. Then
$$
\cC_2(A)=2\sup_{z\in\bc_+}\im z\int_0^1 
{\|(\alpha(x)-zI)^{-1}c(x)G(\ph(x),z)\|^2_{\goth S_2}}\,d\mu(x).
$$
\endproclaim

Analyzing this expression we obtain the following

\proclaim{Corollary 2.9}
Condition~\thetag{2.12} implies~\thetag{UTB}. 
\endproclaim

\head
3. Proofs
\endhead

In this section we prove and comment on the results presented in Section~2.

\subheading{3.1. From the resolvent of the operator $A$ to an auxiliary 
Cauchy problem} First of all we reduce calculation of the resolvent of 
$A$ to some Cauchy problem and then we prove its solvability.

\proclaim{Lemma 3.1} Let $h\in L^2(H,\mu)$ and $z\in\bc_+$ be chosen 
in a way 
that there exists a function $f\in\dom(A^*)$ solving the equation 
$(A^*-zI)f=h$ and the operator-valued function
$\Omega$ defined in~\thetag{2.6} is uniformly bounded on $[0,M]$. 
Then there exists an $E$@-valued function
$g$ solving the Cauchy problem
$$
\cases
g(t,z)'=c_*(t)^*\Omega(t,z)[c_*(t)g(t,z)-h_*(t)]\\
g(M,z)=0.
\endcases
\tag{3.1}
$$
It is uniquely determined by $f$\rom, and the function $f$ itself can be 
recovered from $g$ by the formula
$$
f(x,z)=\Big[\alpha(x)-zI\Big]^{-1}\Big[h(x)-c(x)g(\ph(x),z)\Big]. 
\tag{3.2}
$$
\endproclaim

Notice that, actually, $g$ depends also on $h\in L^2(\mu)$. 
We shall sometimes write $g(t,z,h)$, $f(x,z,h)$  to emphasize a special 
choice of this parameter. On the other hand, we shall write simply 
$g(t)$, $f(x)$, when the parameters $z$ an $h$ are fixed.

The key idea of the calculation is to map $L^2(\mu,[0,1])$
onto a subspace of $L^2(m,[0,M])$, consisting of the functions constant 
on the intervals, which are the ``ranges'' of the point masses of $\mu$. 
According to~\thetag{2.1}, integration with respect to $\mu$ will be 
replaced by integration with respect to Lebesgue measure. We shall come hence
to a differential equation for an auxiliary function $g$.

The proof of existence of the solution to the corresponding Cauchy problem 
is postponed to the next subsection.  

The  argument presented here follows the main lines of~\cite{16}, 
subsect.~2.6--2.10. We let $x$ and $s$ vary in the interval $[0,1]$ and
$t$ and $\tau$ vary in the interval $[0,M]$, respectively.

\demo{Proof} Define a function $g$ by the formula
$$
g(t)=-i\int^M_tc_*(\tau)^*f_*(\tau)d\tau.
$$
Since, by~\thetag{2.1}
$$
\int^1_{x-}c(s)^*f(s)d\mu(s)=\int^M_{\ph(x)}c_*(\tau)^*f_*(\tau)d\tau,
$$
the identity
$$
(\alpha(x)-zI)f(x)-ic(x)\int_{x-}^1 c(s)^*f(s)d\mu(s)=h(x),
\tag{3.3}
$$
can be rewritten in the form
$$
(\alpha(x)-zI)f(x)+c(x)g(\ph(x))=h(x),
\tag{3.4}
$$
and this formula immediately implies~\thetag{3.2}. Hence we need only 
to verify that $g$ defined above satisfies the equation~\thetag{3.1}.

Take a point $t\in[0,M]$ and assume that $\mu(\{x\})>0$ for $x=\psi(t)$.
Recall that in this case $t\in[\ph(x-0),\ph(x+0)]$ and for any function
$v$ on $[0,1]$, we have $v_*(\tau)=v_*(t)=v(x)$ for all 
$\tau\in(\ph(x-0),\ph(x+0))$. Therefore,
$$
g(\ph_*(t))=g(t)-i\int_{\ph_*(t)}^tc_*(\tau)^*f_*(\tau)d\tau=
g(t)-i(t-\ph_*(t))c_*(t)^*f_*(t).
$$
Since $\ph_*(t)=t$ whenever $\mu(\{\psi(t)\})=0$, the latter formula
is valid for every $t$. After plugging it in~\thetag{3.4}, we get
$$
(\alpha_*(t)-zI)f_*(t)+c_*(t)\big[g(t)-i(t-\ph_*(t))c_*(t)^*f_*(t)\big]
=h_*(t),
$$
or
$$
i\Big[i\big(\alpha_*(t)-zI\big)+(t-\ph_*(t))c_*(t)c_*(t)^*\Big]f_*(t)
=c_*(t)g(t)-h_*(t).
$$
Applying the operator $\Omega(t)$ to both sides, we get 
$$
if_*(t)=\Omega(t)[c_*(t)g(t)-h_*(t)].
$$
The definition of $g$ implies that it is an absolutely continuous
function and
$$
g'(t)=ic_*(t)^*f_*(t)=c_*(t)^*\Omega(t)[c_*(t)g(t)-h_*(t)].
$$
i.e., $g$ is a solution to the Cauchy problem~\thetag{3.1}.\qed
\enddemo

\subheading{3.2. Some Cauchy problems, existence of solutions and 
their estimates}
Our main tools in this subsection are the simplest kind of Picard 
approximation and the so-called ``fundamental lemma''. All necessary 
background can be found in~\cite1.
 
\proclaim{Lemma 3.2 \rm{(\cite1, ch.\;2)}} Let $u,v$ be positive 
measurable functions on $[0,M]$, $C_1$ be a non-negative constant and let
$$
u(t)\le C_1+\int_t^M u(\tau)v(\tau)d\tau.
$$
Then
$$
u(t)\le C_1\exp \left(\int^M_tv(\tau)d\tau\right).
$$
\endproclaim

The proof of the following lemma repeats the arguments of Theorem~1
from~\cite1, ch.~1, in a bit more general setting.

\proclaim{Lemma 3.3} Let $\Phi$ and $\Psi$ be two operator-valued 
functions on $[0,M]$\rom. Suppose that they are measurable with respect 
to Lebesgue measure and that the scalar functions $\|\Phi\|$ and $\|\Psi\|$ 
are summable. Then the Cauchy problem
$$
\cases
X(t)'=\Phi(t)X(t)-X(t)\Psi(t)\\
X(M)=X_0
\endcases
\tag{3.5}
$$
is solvable\rom, its solution is unique and absolutely continuous with
respect to the operator norm. Moreover\rom, it satisfies the estimate
$$
\|X(t)\|\le\|X_0\|\exp\left(\int_t^M\left(\|\Phi(\tau)\|+
\|\Psi(\tau)\|\right)d\tau\right).
\tag{3.6}
$$
\endproclaim

\demo{Proof}
We briefly describe the method of Picard approximation. Taking $X_0$ to be
the constant function from the boundary condition of~\thetag{3.5}, we put
$$
X_{k+1}(t)=X_0-
\int^M_t\bigl(\Phi(\tau)X_k(\tau)-X_k(\tau)\Psi(\tau)\bigr)d\tau.
\tag{3.7}
$$
Now, we prove by induction the following estimate
$$
\|X_k(t)-X_{k-1}(t)\|\le\frac{\Gamma(t)^k}{k!}\|X_0\|,
$$
where $\Gamma(t)=\int^M_t\left(\|\Phi(\tau)\|+\|\Psi(\tau)\|\right)d\tau$. 
The estimate is obvious for $k=1$. We make a step of induction for $k>1$
$$
\align
\|X_{k+1}(t)-X_k(t)\|
=&\Bigl\|\int^M_t\Bigl(\Phi(\tau)\big[X_k(\tau)-X_{k-1}(\tau)\big]-
\big[X_k(\tau)-X_{k-1}(\tau)\big]\Psi(\tau)\Bigr)d\tau\Bigr\|\\
\le&\int^M_t\big(\|\Phi(\tau)\|+\|\Psi(\tau)\|\big)
\|X_k(\tau)-X_{k-1}(\tau)\|d\tau\\
\le&\frac1{k!}\int^M_t(-\Gamma(\tau)')\Gamma(\tau)^k\|X_0\|d\tau\\
=&\frac{\Gamma(t)^{k+1}}{(k+1)!}\|X_0\|.
\endalign
$$
Hence, the series
$$
X(t)=\sum_{k=0}^\infty (X_{k+1}(t)-X_k(t))+X_0=\lim X_k(t)
$$
is uniformly convergent with respect to the operator norm topology and
satisfies the required estimate~\thetag{3.6}. Passing to the limit in 
the relation~\thetag{3.7}, we see that the function $X$ satisfies the
equation
$$
X(t)=X_0-\int^M_t\bigl(\Phi(\tau)X(\tau)-X(\tau)\Psi(\tau)\bigr)d\tau
\tag{3.8}
$$
and, therefore, it is absolutely continuous and solves the Cauchy 
problem~\thetag{3.5}.

Notice, that estimate~\thetag{3.6} can be easily proved using
the ``fundamental lemma''. Indeed, since~\thetag{3.8} gives
$$
\|X(t)\|\le\|X_0\|+\int^M_t\|X(\tau)\|\big(\|\Phi(\tau)\|+
\|\Psi(\tau)\|\big)d\tau,
$$
it remains to apply Lemma~3.2 to obtain the bound.

The uniqueness of the solution is a simple consequence of 
estimate~\thetag{3.6}. Indeed, suppose $X_1$ and $X_2$ solve Cauchy 
problem~\thetag{3.5}. Then their difference $X=X_2-X_1$ satisfies the 
same equation~\thetag{3.5} with the boundary condition $X_0=0$. 
Estimate~\thetag{3.6} yields $X(t)=0$ for all $t\in[0,M]$.
\qed 
\enddemo

\proclaim{Corollary~3.4} If the operator $X_0$ is invertible, then the
solution to~\thetag{3.5} is invertible for all $t\in[0,M]$ and its
inverse $Y(t)=X(t)^{-1}$ is the solution to the following Cauchy problem
$$
\cases
Y(t)'=\Psi(t)Y(t)-Y(t)\Phi(t)\\
Y(M)=X_0^{-1}.
\endcases
\tag{3.9}
$$
\endproclaim
   
\demo{Proof}
We take the function $Y$ solving the Cauchy problem~\thetag{3.9} (its values
$Y(t)$ are bounded operators for all $t$ by Lemma~3.3) and check that
$X(t)Y(t)=Y(t)X(t)=I$ for all $t\in[0,M]$. Observe that the function 
$Z=XY$ solves the Cauchy problem
$$
\cases
Z(t)'=\Phi(t)Z(t)-Z(t)\Phi(t)\\
Z(M)=I,
\endcases
$$
and the function $Z(t)=I$ solves the same problem as well. By
Lemma~3.3, the Cauchy problem has unique solution, therefore $X(t)Y(t)=I$.
Replacing $\Phi$ by $\Psi$, we get the equality $Y(t)X(t)=I$ in a similar 
way.\qed
\enddemo

Now, we return to Cauchy problems~\thetag{2.5} and~\thetag{3.1}. 
Before proceed with the proofs, we would like to make a 
remark on existence of the inverse operator $\Omega(t,z)$. 
This operator exists for $\im z\ge\hlf\sup\mu_x\|k(x,x)\|$.
But for the sake of simplicity we shall assume condition~\thetag{2.7} to be
fulfilled, i.e.,
$$
\im z\ge1+\hlf\sup\mu_x\|k(x,x)\|.
$$
Then, $\Omega(t,z)$ are contractions for all $t\in[0,M]$. We state this fact
as a separate lemma.

\proclaim{Lemma 3.5}  We have the
following estimates for $z$ satisfying~\thetag{2.7}\rom:
\roster
\item"i)" $\|\Omega(t,z)\|\le1$;
\item"ii)" $\|c_*(t)^*\Omega(t,z)c_*(t)\|\le\|k_*(t,t)\|$;
\item"iii)" $\|c_*(t)^*\Omega(t,z)c_*(t)\|_{\goth S_1}\le
\|k_*(t,t)\|_{\goth S_1}$.
\endroster
\endproclaim

\demo{Proof}
Denoted by $R(x,z)$ the resolvent of $\alpha(x)$, i.e.,
$R(x,z)=(\alpha(x)-zI)^{-1}$. Then for every $z$ satisfying~\thetag{2.7}
the inverse operator $\Omega(t,z)$ can be defined by any of the following
two expressions:
$$
\align
\Omega(t,z)&=-iR_*(t,z)\big[I-i(t-\ph_*(t))k_*(t,t)R_*(t,z)\big]^{-1}\\
&=-i\big[I-i(t-\ph_*(t))R_*(t,z)k_*(t,t)\big]^{-1}R_*(t,z).
\tag{3.10}
\endalign
$$
Since $\displaystyle\|R_*(t,z)\|\le\frac1{\im z}$ and 
$|t-\ph_*(t)|\le\hlf\mu_x$,
estimating any of these two expressions we get
$$
\|\Omega(t,z)\|\le\frac{\frac1{\im z}}{1-\hlf\mu_x\|k(x,x)\|\frac1{\im z}}
\le1
$$
Two other estimates easily follow from the first one if we take into
account that $\|c(x)\|^2=\|k(x,x)\|$ and $\|c(x)\|^2_{\goth S_2}=
\|k(x,x)\|_{\goth S_1}$:
$$
\align
\|c_*(t)^*\Omega(t,z)c_*(t)\|&=\|c_*(t)^*\|\,\|\Omega(t,z)\|\,\|c_*(t)\|
\le\|k_*(t,t)\|;\\
\|c_*(t)^*\Omega(t,z)c_*(t)\|_{\goth S_1}&=
\|c_*(t)^*\|_{\goth S_2}\|\Omega(t,z)\|\,\|c_*(t)\|_{\goth S_2}
\le\|k_*(t,t)\|_{\goth S_1} \qed
\endalign
$$
\enddemo

Note that for the generator $c_*(t)^*\Omega(t,z)c_*(t)$ we can use
the following ``more symmetric'' formula
$$
c_*(t)^*\Omega(t,z)c_*(t)=c_*(t)^*R_*(t,z)c_*(t)
\Big[I-i\big(t-\ph_*(t)\big)c_*(t)^*R_*(t,z)c_*(t)\Big]^{-1}.
\tag{3.11}
$$

Now, we reformulate Lemma~3.3 and Corollary~3.4 for Cauchy
problem~\thetag{2.5}.

\proclaim{Corollary~3.6} Under assumption~\thetag{2.7}
Cauchy problem~\thetag{2.5} is solvable, its solution is unique and it
satisfies the following bound
$$
\|G(t,z)\|\le\exp\left(\int_t^M\|k_*(\tau,\tau)\|d\tau\right).
$$
Furthermore, $G(t,z)$ is invertible for all $t\in[0,M]$, its inverse
is the unique solution of the Cauchy problem
$$
\cases
\big(G(t,z)^{-1}\big)'&=\;\;-G(t,z)^{-1}c_*(t)^*\Omega(t,z)c_*(t)\\
G(M,z)^{-1}&=\;\;I.
\endcases
$$
The inverse satisfies the same estimate
$$
\|G(t,z)^{-1}\|\le\exp\left(\int_t^M\|k_*(\tau,\tau)\|d\tau\right).
$$
\endproclaim

\demo{Proof} We need only to mention that the 
operator-valued function $\Phi(t)=c_*(t)^*\Omega(t,z)c_*(t)$ satisfies
the estimate $\|\Phi(t)\|\le\|k_*(t,t)\|$ proved in Lemma~3.5. Thus, 
everything except the last claim follows by Lemma~3.3. 
The latter assertion is contained in Corollary~3.4.
\qed
\enddemo

The trace class estimate in Lemma~3.5 allows us to make the following
conclusion.

\proclaim{Corollary~3.7} Under assumption~\thetag{2.7}
the solution of Cauchy problem~\thetag{2.5} has the following 
additional properties\rom: $I-G(t,z)\in\goth S_1$ and
$$
\|I-G(t,z)\|_{\goth S_1}\le
\exp\left(\int_t^M\|k_*(\tau,\tau)\|_{\goth S_1}d\tau\right)-1.
$$
\endproclaim

\demo{Proof}
Rewrite the equation for $G$
$$
G(t,z)=I-\int_t^Mc_*(\tau)^*\Omega(\tau,z)c_*(\tau)G(\tau,z)d\tau
$$
in the form
$$
I-G(t,z)=\int_t^Mc_*(\tau)^*\Omega(\tau,z)c_*(\tau)d\tau-
\int_t^Mc_*(\tau)^*\Omega(\tau,z)c_*(\tau)\big[I-G(\tau,z)\big]d\tau.
$$
We immediately obtain the inequality
$$
\|I-G(t,z)\|_{\goth S_1}\le
\int_t^M\|c_*^*(\tau)\Omega(\tau,z)c_*(\tau)\|_{\goth S_1}
\big(1+\|I-G(\tau,z)\|_{\goth S_1}\big)d\tau.
$$
Using bound iii) from Lemma~3.5, we get
$$
1+\|I-G(t,z)\|_{\goth S_1}\le1+
\int_t^M\|k_*(\tau,\tau)\|_{\goth S_1}
\big(1+\|I-G(\tau,z)\|_{\goth S_1}\big)d\tau.
$$
It remains to apply Lemma~3.2 to complete the proof. \qed
\enddemo

We prove solvability of the vector-valued Cauchy problem~\thetag{3.1} 
now.

\proclaim{Lemma 3.8} Let~\thetag{2.7} be fulfilled and
$h\in L^2(H,\mu)$. Then there exists a unique solution $g$ to Cauchy 
problem~\thetag{3.1}. It can be found by the formula
$$
g(t,z,h)=G(t,z)\int^M_t G(\tau,z)^{-1}c_*(\tau)^*\Omega(\tau,z)
h_*(\tau)d\tau.
\tag{3.12}
$$
The solution is absolutely continuous and uniformly 
bounded on $[0,M]$\rom: 
$$
\|g(t,z,h)\|_{E}\le C_2\|h\|_{L^2(H,\mu)},
$$
where $C_2=\hlf \exp\left(\int^1_0\|k(x,x)\|d\mu(x)\right)
\left[\exp\left(2\int^1_0\|k(x,x)\|d\mu(x)\right)-1\right]$.
\endproclaim

\demo{Proof}
Formula~{3.12} for the solution to Cauchy problem~\thetag{3.1}
is well known (see~\cite1, ch.~1, Theorem~3). The absolute 
continuity of this function (which follows from that one of $G$) allows
us to check that this function solves 
problem~\thetag{3.1} by direct calculation . Uniqueness of $g$ follows from Lemma~3.2.
Indeed, if $g_1$ and $g_2$ are any two solution of~\thetag{3.1}, then 
their difference $g_0=g_2-g_1$ solves the Cauchy problem
$$
g_0(t,z)'=c_*(t)^*\Omega(t,z)c_*(t)g_0(t,z),\qquad g_0(M,z)=0,
$$
and, hence,
$$
g_0(t,z)=\int_t^Mc_*(\tau)^*\Omega(\tau,z)c_*(\tau)g_0(\tau,z)d\tau.
$$
So, the estimate
$$
\|g_0(t,z)\|\le\int_t^M\|k_*(\tau,\tau)\|\,\|g_0(\tau,z)\|d\tau
$$
permits to apply Lemma~3.2 with $C_1=0$.

Using the bounds from Corollary~3.6 and from Lemma~3.5, we 
estimate the norm of $g$ (we use the notation $\omega(t)=
\int^M_t\|k_*(\tau,\tau)\|d\tau$ here)
$$
\align
\|g(t)\|_E&\le\|G(t)\|\int^M_t\|G(\tau)^{-1}\|
\|c_*(\tau)^*\Omega(\tau)h_*(\tau)\|_E d\tau\\
&\le e^{\omega(t)}\int_t^M e^{\omega(\tau)}
\|k_*(\tau,\tau)\|^\hlf\|h_*(\tau)\|_H d\tau\\
&\le e^{\omega(t)}\int_t^M e^{2\omega(\tau)}\big(-\omega(\tau)'\big)d\tau
\|h\|_{L^2(H,\mu)}\\
&=\hlf e^{\omega(t)}\left(e^{2\omega(t)}-1\right)\|h\|\le C_2\|h\|,\qed
\endalign
$$
\enddemo

\subheading{3.3. From the solution of the Cauchy problem to the resolvent
of the operator $A^*$}
The following lemma, together with the previous one, is converse, in 
some sense, to Lemma~3.1. 
In what follows, we shall denote the resolvent of
$\alpha(x)$ again by $R(x,z)$, i.e., $R(x,z)=(\alpha(x)-zI)^\inv$.

\proclaim{Lemma 3.9} The half-plane 
$\im z\ge1+\hlf\sup\big\{\mu_x\|k(x,x)\|\big\}$
does not intersect the spectrum of the operator $A^*$. If
$h\in L^2(H,\mu)$\rom, and $g$ is the solution to Cauchy 
problem~\thetag{3.1}\rom, then the function $f$ defined by~\thetag{3.2} 
is precisely the resolvent of $A^*$ applied to the vector $h$\rom, 
i.e.\rom, $f=(A^*-zI)^\inv h$.
\endproclaim

\demo{Proof}
We start with proving that $f$ defined by~\thetag{3.2} belongs to
$\dom(A^*)$. Taking into account that $\|c(x)^*c(x)\|=\|k(x,x)\|$, we get
$$
\|h(x)-c(x)g(\ph(x),z)\|_H
\le\|h(x)\|_H+C_2\|k(x,x)\|^\hlf\|h\|\in L^2(\mu).
$$
Therefore, the inequality $\displaystyle\|R(x,z)\|\le\frac1{\im z}\le1$ 
implies that $f\in L^2(H,\mu)$. Since $R(x,z)$ maps any vector into 
$\dom(\alpha(x))$ the vector-valued function $\alpha f$ is well defined
and belongs to $L^2(H,\mu)$, because $\alpha(x)R(x,z)=I+zR(x,z)$ and
$R(x,z)$ is uniformly bounded in $x$. Thus, we have proved that 
$f\in\dom(A^*)$.

Now, we verify the relation $(A^*-zI)f=h$.
$$
\align
\big((A^*&-zI)f\big)(x)-h(x)
=-c(x)g(\ph(x))-ic(x)\int_{\ph(x)}^M c_*(t)^*f_*(t)dt\\
&=c(x)\int_{\ph(x)}^M c_*(t)^*\Big(\Omega(t)[c_*(t)g(t)-h_*(t)]
-if_*(t)\Big)dt\\
&=c(x)\int_{\ph(x)}^M c_*(t)^*\Omega(t)\Big(c_*(t)g(t)-h_*(t)\\
&\hskip1cm
-i\big(i(\alpha_*(t)-zI)+(t-\ph_*(t))c_*(t)c_*(t)^*\big)f_*(t)\Big)dt\\
&=c(x)\int_{\ph(x)}^M c_*(t)^*\Omega(t)c_*(t)\Big[g(t)-g(\ph_*(t))-
i(t-\ph_*(t))c_*(t)^*f_*(t)\Big]dt.
\tag{3.13}
\endalign
$$
We see that the expression in the brackets is equal to zero at every $t$
with $\mu(\{\psi(t)\})=0$, because  at these points 
$\ph_*(t)=t$. So, from now on, we may assume that $\mu$ has a mass at the
point $x=\psi(t)$. Recall that in this case $t\in\big[\ph(x-0),\ph(x+0)\big]=
\big[\mu[0,x),\mu[0,x]\big]$. We shall check that  the expression in the 
brackets in~\thetag{3.13} equals zero on the open intervals. The set of 
endpoints of these intervals is at the most countable, i.e., of Lebesgue 
measure zero and may be disregarded.

We want to know more on the behavior of $g$ on the interval 
$(\ph(x-0),\ph(x+0))$.
In other words, we are led to study the behavior of $G$. The 
information we 
obtain will be also useful for 
study of factorizations of the characteristic function $S_A$. 

\proclaim{Sublemma 3.10} If $x$ is a point of non zero mass for
$\mu$, then  the solution $G$ 
of~\thetag{2.5} has the form
$$
G(t,z)=\Big[I-i\big(t-\ph(x)\big)c(x)^*R(x,z)c(x)\Big]G(\ph(x),z).
\tag{3.14}
$$
at the points $t\in[\ph(x-0),\ph(x+0)]$.
\endproclaim

\demo{Proof}. We know that the solution of~\thetag{2.5} is unique
and our expression coincides with the solution at least at the point 
$t=\ph(x)$. So, it is enough to prove that this expression satisfies the
differential equation in~\thetag{2.5}.

The left-hand side of this equation is
$$
G(t,z)'=-i\, c(x)^*R(x,z)c(x)G(\ph(x),z)
$$
and, since $v_*(t)=v(x)$ for any function $v$ and for any interior point 
$t$ of the considered interval, the right-hand side of~\thetag{2.5} is
$$
\align
c_*(t)^*\Omega(t,z)c_*(t)G(t,z)
&=c(x)^*\Omega(t,z)\Big[I-i\big(t-\ph(x)\big)c(x)c(x)^*R(x,z)\Big]
c(x)G(\ph(x),z)\\
&=-i\,c(x)^*R(x,z)c(x)G(\ph(x),z).
\endalign
$$
Thus, $G$ satisfies the equation in~\thetag{2.5} and the sublemma is 
proved. \qed
\enddemo

Now, we return back to the function $g$ on the interval $(\ph(x-0),\ph(x+0))$.
Denote by $I(a,b)$ the integral 
$\int^b_aG(\tau)^{-1}c_*(\tau)^*\Omega(\tau)h_*(\tau)d\tau$. 
Then, by~\thetag{3.12}, we have
$$
\align
g(t)-g(\ph(x))&
=G(t)I(t,M)-G(\ph(x))I(\ph(x),M)\\
&=G(t)I(t,\ph(x))+[G(t)-G(\ph(x))]I(\ph(x),M)\\
&=G(t)I(t,\ph(x))-i\big(t-\ph(x)\big)c(x)^*R(x,z)c(x)
G(\ph(x))I(\ph(x),M)\\
&=G(t)I(t,\ph(x))-i\big(t-\ph(x)\big)c(x)^*R(x,z)c(x)g(\ph(x)).
\tag{3.15}
\endalign
$$
The first summand in~\thetag{3.15} is equal 
$\displaystyle i\big(t-\ph(x)\big)c(x)^*R(x,z)h(x)$. We use 
representation~\thetag{3.10} to show this:
$$
\align
G(t)&I(t,\ph(x))=\Big[I-i\big(t-\ph(x)\big)c(x)^*R(x)c(x)\Big]
G(\ph(x))\times\\
&\hskip1cm\times\int_t^{\ph(x)}G(\ph(x))^{-1}
\Big[I-i\big(\tau-\ph(x)\big)c(x)^*R(x)c(x)\Big]^{-1}c(x)^*
\Omega(\tau)h(x)\,d\tau\\
&=-ic(x)^*\Big[I-i\big(t-\ph(x)\big)R(x)c(x)c(x)^*\Big]\times\\
&\hskip1cm\times
\int_t^{\ph(x)}\Big[I-i\big(\tau-\ph(x)\big)R(x)c(x)c(x)^*\Big]^{-2}
R(x)h(x)\,d\tau\\
&=i\big(t-\ph(x)\big)c(x)^*R(x)h(x).
\endalign
$$
In the latter identity we used the following relation
$$
\int_t^a[I-(\tau-a)X]^{-2}d\tau=-(t-a)[I-(t-a)X]^{-1}
$$
being valid if $\displaystyle \|X\|\le\frac1{|a-t|}$. This condition
is fulfilled for $a=\ph(x)$ and $X=iR(x)k(x,x)$ under 
assumption~\thetag{2.7}.

Substituting this formula in~\thetag{3.15}, we obtain 
$$
\align
g(t)-g(\ph(x))&=i\big(t-\ph(x)\big)c(x)^*R(x)h(x)
-i\big(t-\ph(x)\big)c(x)^*R(x)c(x)g(\ph(x))\\
&=i\big(t-\ph(x)\big)c(x)^*R(x)\big[h(x)-c(x)g(\ph(x))\big]
=i\big(t-\ph(x)\big)c(x)^*f(x).
\endalign
$$
Hence, \thetag{3.13} implies $(A^*-zI)f=h$, i.e., $\Range(A^*-zI)$ is
the whole space $L^2(H,\mu)$. To prove that the resolvent 
$(A^*-zI)^{-1}$ does exist,
we have only to show that $\ker(A^*-zI)=0$. This is
true by Lemma~3.1, because $f$ is uniquely recovered by means of~\thetag{3.2}
from the solution $g$ to Cauchy problem~\thetag{3.1}. The latter 
solution is unique, and the lemma is proved.
\qed
\enddemo

\proclaim{Corollary~3.11}
The operators $A^*$ and $A$ are maximal.
\endproclaim

\demo{Proof}. The dissipative operators $A$ and $-A^*$ are maximal or
not simultaneously. We just proved that $\Range(A^*-zI)=L^2(H,\mu)$
for some $z\in\bc_+$, therefore, $-A^*$ has no dissipative extensions,
i.e., it is maximal. \qed
\enddemo

\subheading{3.4. Proof of Theorem~2.1 and Corollary~2.2}

\demo{Proof of Theorem 2.1}
As in the previous, we shall use the following notation:
$f(\;\cdot\;,z,h)=(A^*-zI)^{-1}h$, the $E$@-valued function $g$ is determined
by $f$ according to the formula
$$
g(t,z,h)=-i\int_t^Mc_*(\tau)^*f_*(\tau,z,h)d\tau.
$$
We remind that  $g$ can also be expressed in terms of the solution $G$ of Cauchy 
problem~\thetag{2.5} with the help of~\thetag{3.12}.

We take an arbitrary vector $e\in E$ and calculate $S_A$ using 
formula~\thetag{2.4}:
$$
\align
S_A(z)e&=e+ic^*(A^*-zI)^{-1}ce=e+ic^*f(\;\cdot\;,z,ce)=
e+i\int_0^1c(x)^*f(x,z,ce)d\mu(x)\\
&=e+i\int_0^Mc_*(t)^*f_*(t,z,ce)dt=e-g(0,z,ce).
\endalign
$$
By Corollary~3.4, the inverse operator $G^{-1}$ satisfies the 
equation $\left(G(t)^{-1}\right)'=-G(t)^{-1}$ $c_*(t)^*\Omega(t)c_*(t)$.
Thus, by~\thetag{3.12}, we have
$$
\align
g(t,z,ce)&=
G(t,z)\int_t^MG(\tau,z)^{-1}c_*(\tau)^*\Omega(\tau,z)c_*(\tau)e\,d\tau\\
&=-G(t,z)\int_t^M\left(G(\tau,z)^{-1}\right)'e\,d\tau\\
&=G(t,z)\left[G(t,z)^{-1}-I\right]e=e-G(t,z)e.
\tag{3.16}
\endalign
$$
Hence, $S_A(z)e=G(0,z)e$ for arbitrary $e$ from $E$. The theorem is proved.
\qed
\enddemo

It is well known that the solutions to the system of differential 
equations~\thetag{2.5} are given by the so called multiplicative integral 
$$
G(t)=\int_t^{M\atop\car} e^{-c_*(\tau)^*\Omega(\tau)c_*(\tau)d\tau},
$$
(see appendix in \cite9 for a detailed discussion of the concept). 
Furthermore, it is possible to calculate $\det G(t)$ as
$$
\det G(t)=\exp
\left\{-\int_t^M\tr\big[c_*(\tau)^*\Omega(\tau)c_*(\tau)\big]d\tau\right\}
$$
(the matrix-valued case is described by Theorem~2 in~\cite1, generalization
to the trace class operators can be found in~\cite3, ch.~IV).

\demo{Proof of Corollary~2.2}
First, we assume~\thetag{2.7} to be fulfilled. We have 
$S_A(z)=G(0,z)$ and, by Corollary~3.7, $I-G(t,z)\in\goth S_1$, so the 
determinant of $S_A(z)$ is well defined
$$
\det S_A(z)=
\exp\left\{-\int_0^M\tr\big[c_*(t)^*\Omega(t,z)c_*(t)\big]dt\right\}.
$$
For every $x$ with $\mu_x>0$, we put 
$E_1=\cup_{x:\,\mu_x>0}[\ph(x-0),\ph(x+0)]$ and $E_2=[0,M]\bsl E_1$.
Recall that $\ph_*(t)=\ph(x)$ for $t\in(\ph(x-0),\ph(x+0))$ and $\ph_*(t)=t$
for $t\in E_2$.  

We take the orthonormal family $\{e_j(x)\}$ in $H$
(see~\thetag{2.3}) diagonalizing $\alpha(x)$ and $k(x,x)$.
Then the operator $\Omega(t,z)$ can be written in the form
$$
\Omega(t,z)=\sum_j\frac1{i(\alpha_{j*}(t)-z)+(t-\ph_*(t))\kappa_{j*}(t)^2}
\big(\;\cdot\;,e_{j*}(t)\big)e_{j*}(t).
$$
Now we consider the family $q_j(x)=\frac1{\kappa_j(x)}c(x)^*e_j(x)$.  
Since $k(x,x)=c(x)c(x)^*$, this is an orthonormal family in $E$ and
$$
c_*(t)^*\Omega(t,z)c_*(t)=\sum_j\frac{\kappa_{j*}(t)^2}
{i(\alpha_{j*}(t)-z)+(t-\ph_*(t))\kappa_{j*}(t)^2}
\big(\;\cdot\;,q_{j*}(t)\big)q_{j*}(t),
$$
whence
$$
\tr c_*(t)^*\Omega(t,z)c_*(t)=\sum_j\frac{\kappa_{j*}(t)^2}
{i(\alpha_{j*}(t)-z)+(t-\ph_*(t))\kappa_{j*}(t)^2}.
$$

First, we consider the integral over $E_1$. Fix a point $x$ with 
$\mu_x>0$. Then
$$
\gather
-\int^{\ph(x+0)}_{\ph(x-0)}\tr\big[c_*(t)^*\Omega(t,z)c_*(t)\big]dt
=-\sum_j\int^{\ph(x+0)}_{\ph(x-0)}\frac{\kappa_j(x)^2\,dt}
{i(\alpha_j(x)-z)+(t-\ph(x))\kappa_j(x)^2}\\
=-\sum_j\int^{\hlf\mu_x}_{-\hlf\mu_x}\frac{d\tau}
{\tau+i\frac{\alpha_j(x)-z}{\kappa_j(x)^2}}
=\sum_j\log\frac{z-z_j(x)}{z-\ovl z_j(x)}.
\endgather
$$
Thus, we have got the Blaschke product in formula~\thetag{2.8},
the unimodular factors $e^{i\phi_j(x)}$ are introduced to make the 
product converge.

Integration over $E_2$ corresponds to integration with respect to the 
continuous part $\mu_c$ of $\mu$. Indeed, let $v\in L^1(\mu)$ and 
$$
\tilde v(x)=\cases
v(x)\quad&\text{if}\ \mu(\{x\})=0,\\
0& \text{otherwise.}
\endcases
$$
Then
$$
\int_0^1v(x)d\mu_c(x)=\int_0^1\tilde v(x)d\mu(x)=\int_0^M\tilde v_*(t)dt
=\int_{E_2}\tilde v_*(t)dt=\int_{E_2}v_*(t)dt.
$$

Since $\Omega(t,z)=-iR(z,\psi(t))$ for $t\in E_2$, we can easily 
calculate the integral over $E_2$:
$$
\align
-\int_{E_2}\tr\big[c_*(t)^*\Omega(t,z)c_*(t)\big]dt
&=i\,\int_{E_2}\tr\big[c_*(t)^*R(\psi(t),z)c_*(t)\big]dt\\
&=i\int_0^1\tr\big[c(x)^{*}R(x,z)c(x)\big] d\mu_c(x).
\endalign
$$
So,  we proved formula~\thetag{2.8} for $\im z\ge1+\hlf\sup\mu_x\|k(x,x)\|$.
Since both sides of \thetag{2.8} are $H^\infty$@-functions on $\bc_+$,
analytic continuation finishes the proof.\qed 
\enddemo

\subheading{3.5. Proof of Theorem~2.3}
The proof of Theorem~2.3 is based on the following lemma.

\proclaim{Lemma 3.12}
Condition~\thetag{2.10} of Theorem~{\rm2.3} \;is necessary, and under
assumption~\thetag{2.2} sufficient, for the estimate
$$
\inf_{z\in\bc_+}\{\det S_A(z)\}_{sing,\,out}\ge0.
$$
\endproclaim

\demo{Proof}
First we assume that condition~\thetag{2.2} is fulfilled.
Formula~\thetag{2.8} shows that the product of the inner singular and 
the outer factors  $\{\det S_A\}_{sing,\, out}$ has the form
$$
\{\det S_A(z)\}_{sing,\,out}=
\exp\Big(i\int_0^1\tr\big[c(x)^*R(x,z)c(x)\big] d\mu_c(x)\Big).
$$ 
Using spectral representation of $\alpha(x)$ we obtain
$$
\{\det S_A(z)\}_{sing,\,out}=
\exp\Big(i\int_0^1\int_\br\frac{\tr[c(x,\lambda)^*c(x,\lambda)]}{\lambda-z}
d\rho_x(\lambda)\,d\mu_c(x)\Big).
$$ 
Since $\tr[c(x,\lambda)^*c(x,\lambda)]\ge0$ and
$
\displaystyle
\re\frac{i}{\lambda-z}=-\frac{\im z}{|\lambda-z|^2},
$
the following identity holds
$$
\align
|\{\det S_A(z)\}_{sing,\; out}|^{-1}&=
\exp\Big(\im z\int_0^1\int_\br\frac{\tr[c(x,\lambda)^*c(x,\lambda)]}
{|\lambda-z|^2}d\rho_x(\lambda)\,d\mu_c(x)\Big)\\
&=\exp\Big(\im z\int_\br\frac{d\nu_c(s)}{|s-z|^2}\Big).
\endalign
$$
The claim of the lemma follows now from the Fatou's theorem on the boundary 
values of the Poisson integral.

If condition~\thetag{2.2} does not fulfilled, then the integral over $E_1$
in the proof of Corollary~2.2 could be not a pure Blaschke product, but
it could have some outer and singular inner factor. Nevertheless, these
additional factors have the modulus not greater than one. Therefore,
in this situation the boundedness of the above Poisson integral remains
necessary for the function $\{\det S_A(z)\}_{sing,\,out}$ to be bounded
away from zero.
\qed
\enddemo

Note that under assumption~\thetag{2.2} we can give another expression
for the measure $\nu_c$:
$$
\nu_c(F)=\sum_j\int_{\alpha^\inv(F)}\kappa_j(x)^2 d\mu_c(x)
$$

Theorem 2.3 immediately follows from Theorem~1.3 and the lemma
we just proved. Indeed, if~\thetag{UTB} is fulfilled, then by Theorem~1.3
the condition~\thetag{1.4} is necessary for the operator $A$ to be similar
to a normal operator, and by Lemma~3.12 for the latter estimate 
condition~\thetag{2.10} is necessary as well.

If the measure $\mu$ is continuous, then Lemma~3.12 yields
estimate~\thetag{1.4} and Theorem 1.3 guarantees required similarity.
\qed

\subheading{3.6. Proof of Theorem~2.4}

\proclaim{Lemma 3.13}
Let condition~\thetag{2.2} be fulfilled. Then inequality~\thetag{2.11} 
holds if and only if the measure
$$
\sigma=\sum_{j,\,x:\,\mu(\{x\})>0}\im z_j(x)\;\delta_{z_j(x)}
\tag{3.17}
$$
is a Carleson one.
\endproclaim

\demo{Proof}
For an arbitrary square $Q=[x_0-h,x_0+h]\times i[0,2h]$ with 
$x_0\in\br$, $h>0$, we have 
$$
\sigma(Q)=\sum_{z_j(x)\in Q}\im z_j(x)=
\hlf\sum_
{\scriptstyle x_0-h\le\alpha_j(x)\le x_0+h,\atop
\scriptstyle 0\le\hlf\mu_x\kappa_j(x)^2\le 2h}
\mu_x\kappa_j(x)^2
=\nu_{d,h}([x_0-h,x_0+h]).
$$
Therefore, the Carleson condition for the measure $\sigma$ has the form
$$
\nu_{d,h}([x_0-h,x_0+h])\le Ch,
$$ 
Thus, condition~\thetag{2.11} is none other than the assertion that
this is a Carleson measure. \qed
\enddemo

Since the (LRG) condition is necessary for similarity to a normal operator,
the above lemma with Theorems~1.6 proves Theorem~2.4.\qed

\subheading{3.7. Proof of Theorem 2.5}
We verify the  conditions of Theorem ~1.3.
Observe firstly that condition i) of Theorem~1.3 is equivalent
to~\thetag{2.10} by Lemma~3.12, and  it is fulfilled under our assumption. 
Hence, we only need to check that the eigen spaces
form an unconditional basis in their span.

By Lemma~3.13, condition~\thetag{2.11} implies that the set of eigenvalues 
$\{z_j(x)\}$ is an $N$-Carleson set. Since $\Lambda=\sigma_p(A)\cap\bc_+$ 
is sparse, it is a Carleson set (if we do not count multiplicities), i.e., 
the Blaschke product with simple zeros from $\Lambda$ satisfies the Carleson 
condition. Since $A$ has no root vectors, this Blaschke product is the 
minimal function of the restriction of $A$ to 
$\spn\{\ker (A-zI):\,z\in\sigma_p(A)\cap\bc_+\}$,
and the eigen spaces form an unconditional basis in their span.
\qed

\subheading{3.8. Proof of Theorem 2.6}
It is convenient to split condition~\thetag{2.12} in
\thetag{2.10} and \thetag{2.11}. Lemma~3.12 implies that condition~i)
of Theorem~1.3 is equivalent to~\thetag{3.10}. By Lemma~3.13,
condition~\thetag{2.11} means that \thetag{3.17} is a Carleson measure.
This is equivalent, together with the sparseness property of $\{z(x)\}$,  
to the Carleson property of the mentioned sequence (see Remark~1.5).
But for an operator with a scalar characteristic function this is
equivalent to the assertion that the operator has no root vectors and
its eigenvectors form an unconditional basis in their span, i.e., is
equivalent to condition~ii) of Theorem~1.3.
The theorem is proved.\qed

\subheading{3.9. Proof of Lemma~2.7 and pertinent observations}
We precede the discussions of factorizations of the characteristic 
function $S_A$ with the following essential refinement of Corollary~3.6.

\proclaim{Lemma 3.14} The solution $G(t,z)$ to the Cauchy
problem~\thetag{2.5} has analytic continuation in parameter $z$ onto 
the whole half-plane $\bc_+$. Moreover\rom, its 
values are contractions for all $t\in[0,M]$\rom, $\im z>0$.
For $t=\ph(x\pm0)$\rom, $x\in[0,1]$\rom, the operator-valued functions 
$G(0,z)G(t,z)^{-1}$ also have analytic continuation to the whole 
half-plane $\bc_+$\rom, their values are contractions\rom,
and the factorization 
$$
G(0,z)=\Big[G(0,z)G(t,z)^{-1}\Big]G(t,z)
$$
is regular {\rm(see, e.g., \cite{12} for the notion of regularity).}
\endproclaim

\demo{Proof} Fix a point $x_0\in(0,1)$ and put $t_0=\ph(x_0+0)=\mu([0,x_0])$.
Let $\mu_0$ be the restriction of the measure $\mu$ to the interval 
$(x_0,1]$, i.e., $\mu_0(F)=\mu(F\cap(x_0,1])$ for any measurable subset 
$F$ of $[0,1]$. We consider the operator $A_0$ in $L^2(\mu_0)$ defined 
by the same formula~\thetag{0.1} with $\mu$ is replaced by $\mu_0$.
It is clear that corresponding Cauchy problem~$\thetag{2.5}_0$ on
the interval $[0,M_0]$, $M_0=M-t_0$, is simply the original Cauchy problem 
on the interval $[t_0,M]$, shifted to the left by $t_0$. 
Therefore, $S_{A_0}(z)=G(t_0,z)$, being the characteristic function of 
a maximal dissipative operator, is analytic and contractive-valued in 
the whole half-plane $\bc_+$. The same conclusion holds for 
$t_0=\ph(x_0-0)=\mu([0,x_0))$, we only have to consider $\mu_0$ defined as 
the restriction of $\mu$ to the closed interval $[x_0,1]$.

If we identify $L^2(\mu_0)$ with the subspace of $L^2(\mu)$ consisting
of functions vanishing on $[0,x_0]$, then the subspace is invariant for 
$A$, and $A_0$ is just the restriction of $A$ to $L^2(\mu_0)$. Further, 
$S_A(z)=[G(0,z)G(t_0,z)^{-1}]G(t_0,z)$ 
is the corresponding regular factorization of the characteristic function,
i.e., the factor $G(0,z)G(t_0,z)^{-1}$ is an analytic contractive-valued
function of $z\in\bc_+$.

Thus, to complete the proof we need to show that $G(t,z)$ has 
contractive-valued analytic
continuation in $z$ if $t\in\big(\ph(x-0),\ph(x+0)\big)$ for every point 
$x$ with $\mu_x>0$. Now, we fix one of these $x$ and use formula~\thetag{3.14} 
proved in Sublemma~3.10. We have

$$
G(\ph(x+0),z)=\Big[I-i\hlf\mu_xc(x)^*R(x,z)c(x)\Big]G(\ph(x),z),
$$
and therefore
$$
G(t,z)=
\Big[I-i\big(t-\ph(x)\big)c(x)^*R(x,z)c(x)\Big]
\Big[I-i\hlf\mu_xc(x)^*R(x,z)c(x)\Big]^\inv G(\ph(x+0),z).
\tag{3.18}
$$
Since
$$
\align
\re\Big[I-i\hlf\mu_xc(x)^*R(x,z)c(x)\Big]
&=I+\hlf\mu_xc(x)^*[\im R(x,z)]c(x)\\
&=I+\hlf\im z\mu_xc(x)^*R(x,z)R(x,z)^*c(x)\ge I
\endalign
$$
for all $z\in\bc_+$, the inverse operator exists and is analytic in
the whole half-plane $\bc_+$. 

Now, we estimate the norm of $G(t,z)$:
$$
\align
I&-\Big[I+i\hlf\mu_xc(x)^*R(x,z)^*c(x)\Big]^\inv
\Big[I+i\big(t-\ph(x)\big)c(x)^*R(x,z)^*c(x)\Big]\times\\
&\hskip1cm\times\Big[I-i\big(t-\ph(x)\big)c(x)^*R(x,z)c(x)\Big]
\Big[I-i\hlf\mu_xc(x)^*R(x,z)c(x)\Big]^\inv\\
&=\left(\hlf\mu_x-\big(t-\ph(x)\big)\right)
\Big[I+i\hlf\mu_xc(x)^*R(x,z)^*c(x)\Big]^\inv
c(x)^*R(x,z)^*\times\\
&\hskip1cm\times\Big\{2\im zI+
\big[\hlf\mu_x+\big(t-\ph(x)\big)\big]c(x)c(x)^*\Big\}\times\\
&\hskip1cm R(x,z)c(x)\Big[I-i\hlf\mu_xc(x)^*R(x,z)c(x)\Big]^\inv\ge0
\endalign
$$
for all $z\in\bc_+$ and $t\in\big[\mu([0,x)),\mu([0,x])\big]$. Therefore,
$$
G(t,z)^*G(t,z)\le G(\ph(x+0),z)^*G(\ph(x+0),z)\le I,
$$
and the proof is complete.\qed
\enddemo

\demo{Proof of Lemma 2.7} By Lemma~3.14, taking two regular factorizations
of $G(0,z)$ corresponding to the points $t_\pm=\ph(x\pm0)$ we get the
following regular factorization 
$$
S_A(z)=S_{x-}(z)B_x(z)S_{x+}
$$
with $B_x(z)=G(\ph(x-0),z)G(\ph(x+0),z)^\inv$.
From~\thetag{3.18}, we obtain 
$$
B_x(z)=\Big[I+\frac{i}2\mu_xc(x)^*R(x,z)c(x)\Big]
\Big[I-\frac{i}2\mu_xc(x)^*R(x,z)c(x)\Big]^{-1}.
$$
Under assumption~\thetag{2.2} we see that $B_x$ is a diagonal 
Blaschke--Potapov product. Indeed, if $\{e_j(x)\}$ is the family of 
joint normalized eigenvectors of $k(x,x)$ and $\alpha(x)$ corresponding 
to the nonzero eigenvalues $\kappa_j(x)^2$ of $k(x,x)$ (sf.~\thetag{2.3}),
then, as in the proof of the Corollary~2.2, we can introduce the 
orthonormal family $\{q_j(x)\}$ in $E$, 
$q_j(x)=\frac1{\kappa_j(x)}c(x)^*e_j(x)$.
If we denote by $Q(x)$ the orthogonal projection of $E$ onto orthogonal
complement to $\spn_j\{q_j(x)\}$, then
$$
\align
B_x(z)&=\Big[I+\frac{i}2\sum_j\frac{\kappa_j(x)^2}{\alpha_j(x)-z}
\big(\;\cdot\;,q_j(x)\big)q_j(x)\Big]
\Big[I-\frac{i}2\sum_j\frac{\kappa_j(x)^2}{\alpha_j(x)-z}
\big(\;\cdot\;,q_j(x)\big)q_j(x)\Big]\\
&=Q(x)+\sum_j\frac{z-z_j(x)}{z-\ovl{z_j(x)}}
\big(\;\cdot\;,q_j(x)\big)q_j(x).
\endalign
$$
The lemma is proved.\qed
\enddemo

\proclaim{Corollary 3.15} For any finite set $x_j\in[0,1]$\rom, 
$1\le j\le n$\rom, there exists a regular factorization of the 
characteristic function of the form 
$$
S_A(z)=S_{x_1-}B_{x_1}S_{x_1,x_2}B_{x_2}S_{x_2,x_3}\cdot\dots\cdot
B_{x_n}S_{x_n+},
\tag{3.19}
$$
where $S_{a,b}=G(\ph(a+0),z)G(\ph(b-0),z)^{-1}$.
\endproclaim

In the ``degenerate'' case $x_1=0$ (or $x_n=1$) the left factor $S_{0-}$
(or the right factor $S_{1+}$, respectively) is the identity operator.

\demo{Proof} Recall that $G(\ph(x+0),z)$ is the characteristic
function of the operator $A_0$ with the measure $\mu_0$ defined as the
restriction of $\mu$ to $(x,1]$ (see the proof of Lemma~3.14).
Now, we apply Lemma~3.14 recursively: first to the
point $x_1$, then, taking $S_{x_1+}$ instead of $S_A$, to the point $x_2$
and so on. \qed
\enddemo

Note that under condition~\thetag{2.2}
we are able to calculate the determinant of every term 
in~\thetag{3.19}:
$$
\align
\det B_x(z)&=\prod_j 
\left(\frac{z-z_j(x)}{z-\ovl{z_j(x)}}e^{i\phi_j(x)}\right),\\
\det S_{a,b}(z)&=\prod_{j,\,x:a<x<b} 
\left(\frac{z-z_j(x)}{z-\ovl{z_j(x)}}e^{i\phi_j(x)}\right)\cdot
\exp\left(i\int_a^b\tr [c(x)^*(\alpha(x)-z)^\inv c(x)]d\mu_c(x)\right).
\endalign
$$
This helps us, for example, in ``localization'' of the zeros of the
characteristic function. Fix a $\lambda\in\sigma_p(A)\cap\bc_+$. Since 
zeros of $\det S_A$ satisfy the Blaschke condition, there exists a 
finite set, say $\{x_l\}_{l=1}^n$, of the points from $[0,1]$ such that 
$z_j(x_l)=\lambda$ for some $j$. Consider factorization~\thetag{3.19}
with the chosen set $\{x_l\}$. Then every factor $B_{x_l}$ has a non empty 
kernel at the point $\lambda$ of dimension, say $\ka_l$, so, 
$\ka(\lambda)=\sum^n_{l=1}\ka_l$ is the multiplicity of zero of $\det S_A$ 
at $\lambda$. The function $\det S_A^{-1}b_\lambda^{\ka(\lambda)}$ is
analytic in a neighborhood of $\lambda$ and therefore all factors
$S$ in~\thetag{3.19} are invertible there, and 
the inverses are well defined at $\lambda$.
Let $P_\lambda(x)$ be
the orthogonal projection of $E$ onto $\ker B_x(\lambda)=
\ker\Big[I+\frac{i}2\mu_xc(x)^*R(x,\lambda)c(x)\Big]$, then
$\rank P_\lambda(x_l)=\ka_l$, and $B_{x_l}$ admits factorization
$$
B_{x_l}(z)=B_{x_l}^\lambda(z)
\Big[b_\lambda(z) P_\lambda(x_l)+\big(I-P_\lambda(x_l)\big)\Big],
$$
where the factor $B_{x_l}^\lambda(z)$ is already invertible in a 
neighborhood of the point $z=\lambda$. We denote by $S_{x+}^\lambda$
the ``tail'' of the decomposition~\thetag{3.19}, where all $B_{x_l}$ are
replaced by $B_{x_l}^\lambda$:
$$
S_{x_l+}^\lambda=S_{x_l,x_{l+1}}B_{x_{l+1}}^\lambda S_{x_{l+1},x_{l+2}}
\cdot\dots\cdot B_{x_n}^\lambda S_{x_n+}.
$$
Then, $S_{x_l}^\lambda(z)$ are also invertible in a neighborhood of 
the point $z=\lambda$.

Now, we can analyze the geometric properties of the eigen and root spaces. 
The operator $A$ does not have root vectors if and only if $S_A(z)^{-1}$ 
has only simple poles at all points $\lambda\in\sigma_p(A)\cap\bc_+$, or, 
in other words, $\dim\ker S_A(\lambda)=\ka(\lambda)$. It is clear that the 
latter identity occurs if and only if $\Range P_\lambda(x_l)\subset
\Range S_{x_l+}(\lambda)$.

Assuming that $A$ has no root vectors, we can remove from the decomposition
inverse to~\thetag{3.19} all terms with the poles of the orders greater
than one. We obtain the following expression:
$$
S_A(z)^{-1}=\sum_{l=1}^n S_{x_l+}^\lambda(z)^{-1}
\big[b_\lambda(z)^{-1}P_\lambda(x_l)+(I-P_\lambda(x_l))\big]
B^\lambda_{x_l}(z)^{-1}S_{x_l-}^\lambda(z)^{-1}
$$
We get the following expression for the kernel of $S_A(\lambda)$ by taking the
range of the residue of $S_A^{-1}$ at the point:
$$
\ker S_A(\lambda)=\spn\big\{S_{x_l+}^\lambda(z)^{-1}
\Range P_\lambda(x_l):\,1\le l\le n\big\}.
$$
Note that since $\dim\ker S_A(\lambda)=\ka(\lambda)$ and the subspaces 
in this span have dimensions $\ka_l$, they are linearly independent.

\subheading{3.10. Proof of Lemma 2.8}
Since
$$
\align
\cC_2(A)&=\sup_{z\in\bc_+} 4\im z\cdot \tr [(A^*-zI)^{-1}\im A 
(A-\ovl zI)^{-1}]\\
&=\sup_{z\in\bc_+} 2\im z\cdot \tr [(A^*-zI)^{-1}cc^*
(A-\ovl zI)^{-1}]\\
&=\sup_{z\in\bc_+} 2\im z\cdot\|(A^*-zI)^{-1}c\|^2_{\goth S_2},
\endalign
$$
we have only to calculate 
$$
(A^*-zI)^{-1}ce=f(\;\cdot\;,z,ce)=
R(\;\cdot\;,z)\big(ce-cg(\ph(\;\cdot\;),z,ce)\big),
\qquad e\in E,
$$
(we use the notation introduced in Lemma~3.1). It was checked 
in~\thetag{3.16} that $g(t,z,ce)=e-G(t,z)e$. Therefore, 
$ce-cg(\ph,z,ce)=cG(\ph,z)e$, and choosing an arbitrary orthonormal 
basis $\{e_j\}$ in $E$, we get
$$
\align
\|(A^*-zI)^{-1}c\|^2_{\goth S_2}
&=\sum\|(A^*-zI)^{-1}ce_j\|^2_{L^2(H,\mu)}\\
&=\sum\int^1_0\|R(x,z)c(x)G(\ph(x),z)e_j\|_H^2d\mu(x)\\
&=\int^1_0\|R(x,z)c(x)G(\ph(x),z)\|_{\goth S_2}^2 d\mu(x).
\qed
\endalign
$$

\subheading{3.11. Proof of Corollary 2.9}
We can simply omit $G(\ph(x),z)$ in the integral with respect to the 
continuous part of $\mu$, because 
by Lemma~3.14, these operators are contractions:
$$
\align
\im z\int^1_0&\|R(x,z)c(x)G(\ph(x),z)\|_{\goth S_2}^2 d\mu_c(x)
\le\im z\int^1_0\|R(x,z)c(x)\|_{\goth S_2}^2 d\mu_c(x)\\
&\le\im z\int^1_0\int_\br\frac{\|c(x,\lambda)\|_{\goth S_2}^2}
{|\lambda-z|^2}d\rho_x(\lambda)\,d\mu_c(x)
=\im z\int_\br\frac{d\nu_c(\lambda)}{|\lambda-z|^2}.
\endalign
$$
Under assumptions~\thetag{2.12} (and therefore,~\thetag{2.10}) this is 
the Poisson integral of a measure with
a bounded density, hence it is uniformly bounded.

Now, consider the discrete part of the measure. First, we use 
relation~\thetag{3.14} for $t=\ph(x+0)$ and then we estimate 
$\|G(\ph(x+0),z)\|$ by 1:
$$
\align
\im z\sum_{x:\mu_x>0}
&\|R(x,z)c(x)G(\ph(x),z)\|_{\goth S_2}^2\mu_x\\
&\le\im z\sum_{x:\mu_x>0}
\|R(x,z)c(x)G(\ph(x),z)\Big[I-\frac{i}2\mu_xc(x)^*R(x,z)c(x)\Big]
\|_{\goth S_2}^2\mu_x\\
&=\im z\sum_{x:\mu_x>0}\sum_j\frac{\mu_x \kappa_j(x)^2}
{|\alpha_j(x)-z-\frac{i}2\mu_x\kappa_j(x)^2|^2}\\
&=2\sum_{j,\,x:\mu_x>0}\frac{\im z\im z_j(x)}{|z-\ovl z_j(x)|^2}.
\endalign
$$
The latter expression is bounded, because
by Lemma~3.13 the measure~\thetag{3.17} is Carleson and it remains to
refer to Theorem~1.4 to complete the proof.\qed

\subheading{3.11. An example}
At last, we would like to say that the following question is still open:
what is a condition for a dissipative operator having the Linear Resolvent
Growth property and being a trace class perturbation of a self-adjoint
operator to be similar to a normal operator? We underline that the 
(UTB) condition is not necessarily fulfilled under these assumptions. 
Here we give an example of a normal operator given by formula~\thetag{0.1}
without the (UTB) property.

Consider the scalar case $H=E=\bc$.
Take an arbitrary countable subset $\{x_n\}_{n=1}^\infty$ of the unit interval
$(0,1)$ and a summable sequence of point masses $\mu_n=\mu_{x_n}$. Further, 
put $k(x_n,x_m)=w_n\delta_{n,m}$, $w_n>0$, $\sum w_n\mu_n<\infty$, and pick 
an arbitrary function $\alpha$ finite at the points $x_n$. Now, define an 
operator $A$ by formula~\thetag{0.1}.
This operator is normal. The complete orthogonal family of functions
$\delta_n$, $\delta_n(x_m)=\delta_{n,m}$, is the family of its eigenvectors,
and the corresponding eigenvalues are $z_n=\alpha(x_n)+\hlf iw_n\mu_{x_n}$. 
Indeed,
$$
(A\delta_n)(x)=\alpha(x_n)\delta_n(x)+i\int^{x+}_0k(x,s)\delta_n(s)\,d\mu(s)
=z_n\delta_n(x)
$$
Since the characteristic function $S_A$ is diagonal, we see that the (UTB)
condition is reduced to the property of eigenvalues $z_n$ to form a finite
union of Carleson sequences
$$
\sup_{z\in\bc_+}\tr(I-S_A(z)^*S_A(z))=
\sup_{z\in\bc_+}\;4\sum\frac{\im z\im z_n}{|z-\ovl{z_n}|^2}.
$$
Therefore, every Blaschke sequence $\{z_n\}$ that is not a finite union of 
Carleson sequences supplies us with an example of a normal dissipative trace
class perturbation of a self-adjoint operator without the (UTB) property.

\medskip
\noindent{\bf Acknowledgment.} \ We wish to thank Professor N.~Nikolski 
for valuable discussions on the subject. We are 
grateful to Thomas V.~Pedersen for proof reading the paper.

\end